\documentclass{amsart}

\usepackage{amssymb}
\usepackage{amsfonts}
\usepackage{amsmath}
\usepackage{graphicx}%
\usepackage{lastpage}
\usepackage[utf8]{inputenc}
\usepackage[legalpaper,bookmarks=true,colorlinks=true,linkcolor=blue,citecolor=blue]{hyperref}
\usepackage{fancyhdr}
\usepackage{color}
\usepackage[mathlines]{lineno}
\usepackage{lscape}
\usepackage{epsfig}
\usepackage{natbib}
\usepackage{float}

\newtheorem{theorem}{Theorem}
\theoremstyle{plain}

\newtheorem*{lemma}{Lemma}

\numberwithin{equation}{section}

\begin{document}
\Large
\title[A simple Proof of the Theorem of Sklar]{A simple proof of the theorem of Sklar and its extension to distribution functions}

\author{Gane Samb LO}

\begin{abstract}
In this note we provide a quick proof of the Sklar's Theorem on the existence of copulas by using the generalized inverse functions as in the one dimensional case, but a little more sophisticated.

\bigskip \noindent LERSTAD, Gaston Berger University, Saint-Louis, S\'en\'egal (main affiliation).\newline
\noindent LSTA, Pierre and Marie Curie University, Paris VI, France.\newline
\noindent AUST - African University of Sciences and Technology, Abuja, Nigeria.\newline

\noindent \textit{Corresponding author}. Gane Samb Lo. Email :
gane-samb.lo@edu.ugb.sn, gslo@aust.edu.ng, ganesamblo@ganesamblo.net\\
Permanent address : 1178 Evanston Dr NW T3P 0J9,
Calgary, Alberta, Canada.
\end{abstract}

\maketitle

\section{Introduction} 
\noindent Let us begin by defining generalized functions. Let $[a,b]$ and $[c,d]$ be non-empty intervals of $\mathbb{R}$ and let $G : [a,b] \mapsto [c,d]$ be a non-decreasing mapping such that
\begin{eqnarray*}
c&=& \inf_{x \in [a,b]} G(x),\ \ \ (L11)\\
d&=& \sup_{x \in [a,b]} G(x). \ \ \ (L12)\\
\end{eqnarray*}

\noindent Since $G$ is a mapping, this ensures that
\begin{eqnarray*}a&=&\inf \{x\in \mathbb{R}, \ G\left(x\right) > c \}, \ \ \ (L13)\\
b&=&\sup \{x\in \mathbb{R}, G\left( x\right) < d\}. \ \ \ (L14)\\
\end{eqnarray*}

\noindent If $x=a$ or $x=b$ is infinite, the value of $G$ at that point is meant as a limit. If $[a,b]$ is bounded above or below in $\mathbb{R}$, $G$ is extensible on $\mathbb{R}$ by taking $G\left( x\right)=G\left( a+\right)$ for $x\leq a$ and $G\left( x\right) =G\left( b-0\right) $ for $x\geq b$. As a general rule, we may consider $G$ simply as defined on $\mathbb{R}$. In that case, $a=lep(G)$ and $b=uep(G)$ are called \textit{lower end-point} and \textit{upper end-point} of $G$.\\

\noindent The generalized inverse function of $G$ is given by

$$
\forall u \in [lep(G),uep(G)],  \ G^{-1}\left(u\right) =\inf \left\{ x\in \mathbb{R}, \ G\left( x\right) \geq u\right\}.
$$

\bigskip \noindent The properties of $G^{-1}$ have been thoroughly studied, in particular in \cite{billingsley}, \cite{resnick}. The results we need in this paper are gathered and proved in \cite{ips-wcia-ang} or in \cite{ips-wcia-fr} (Chapter 4, Section 1) and reminded as below.

\begin{lemma} \label{lemA} Let $G$ be a non-decreasing right-continuous function with the notation above. Then  $G^{-1}$ is left-continuous and we have

\begin{equation*}
\forall u \in [c,d], \ G(G^{-1}(u))\geq u \ (A) \  and \  \forall x \in [a,b], \ G^{-1}(G(x))\leq x \ (B) 
\end{equation*}

\bigskip \noindent and 

\begin{equation}
\forall x \in [lep(G), uep(G)], \ G^{-1}(G(x)+0)=x. \label{FF}
\end{equation}
\end{lemma}

\noindent \textit{Proof}. The proof of Formulas (A) and (B) are well-known and can be found in the cited books above. Let us prove Formula \eqref{FF} for any $x \in [a,b]$.\\

\noindent On one side, we start by the remark that $G^{-1}(G(x)+0)$ is the limit of $G^{-1}(G(x)+h)$ as $h \searrow 0$. But for any $h>0$, $G^{-1}(G(x)+h)$ is the infimum of the set of $y\in [a,b]$ such that $G(y)\geq G(x)+h$. Any these $y$ satisfies $y\geq x$. Hence $G^{-1}(G(x)+0)\geq x$.\\

\noindent On the other side $G(x+h) \searrow G(x)$ by right-continuity of $G$, and by the existence of the right-hand limit of the non-decreasing function $G^{-1}(\circ)$, $G^{-1}(G(x+h)) \searrow G^{-1}(G(x)+0)$. Since $G^{-1}(G(x+h))\leq x+h$ by Formula (B), we get that $G^{-1}(G(x)+0)\leq x$ as $h \searrow 0$. The proof is complete. $\square$\\

\bigskip \noindent The remainder of the note will focus on distribution functions (\textit{df}'s) and copulas on $\mathbb{R}^d$. For an introduction to \textit{df}, we refer to \cite{ips-mestuto-ang} (Chapter 11) and for copulas to \cite{nelsen}.\\

\noindent Let us also remind the definition of a distribution function on $\mathbb{R}^d$, $d\geq 1$. A mapping $\mathbb{R}^d \mapsto \mathbb{R}$ is a \textit{df} if and only if : \\

\noindent (\textbf{DF1}) $F$ is right-continuous\\

\noindent and \\

\noindent (\textbf{DF2}) $F$ assigns to non-negative volumes to any cuboid $]a,b]$ (with $a=(a_{1},...,a_{d})\leq b=(b_{1},...,b_{d})$ meaning $a_i\leq b_i$, for all $1\leq i \leq d$), that is :                                                                                              

\begin{equation}
\Delta F(a,b)=\sum_{\varepsilon \in \{0,1\}^d}(-1)^{s(\varepsilon)}F(b+\varepsilon \ast (a-b))\geq 0, \label{VP}
\end{equation}

\bigskip \noindent where 

\begin{equation*}
(x,y)\ast (X,Y)=(x_1X_1, x_2X_2, ..., y_kY_k),
\end{equation*}

\noindent $\varepsilon=(\varepsilon_1,\cdots,\varepsilon_d)$ runs over $\{0,1\}^d$ and $s(\varepsilon)=\varepsilon_1+\cdots+\varepsilon_d$.\\

\bigskip \noindent It becomes a cumulative distribution function \textit{cdf} or simply a probability distribution function if $F$ satisfies the third limit conditions (which by the way ensure that $F$ is non-negative) :\\

\noindent (\textbf{DF3a})
\begin{equation*}
\lim_{\exists i,1\leq i\leq k,t_{i}\rightarrow -\infty }F(t_{1},...,t_{k})=0,
\end{equation*}

\bigskip \noindent and\\

\noindent (\textbf{DF3b})

\begin{equation*}
\lim_{\forall i,1\leq i\leq k,t_{i}\rightarrow +\infty }F(t_{1},...,t_{k})=1.
\end{equation*}

\bigskip \noindent Now we come to copulas. By definition, a copula on $\mathbb{R}^d$ is a \textit{cdf} $C$ whose marginal \textit{cdf}'s defined by, for $1\leq i \leq d$,
$$
\mathbb{R} \ni s \mapsto C_{i}(s)=C\left( +\infty, ..., +\infty,\underset{i-th \ argument}{\underbrace{s}}, +\infty, ..., +\infty\right),
$$
 
\bigskip  \noindent are all equal to the $(0,1)$-uniform \textit{cdf} which in turn is defined by
$$
x \mapsto x 1_{[0,1[} + 1_{[1,+\infty[},
$$

\noindent and we may also write, for all $s \in [0,1]$, 

\begin{equation}
C_{i}(s)=C\left(1, ..., 1,\underset{i-th \ argument}{\underbrace{s}}, 1, ...,1 \right)=s. \label{probcop}
\end{equation}

\bigskip  \noindent Based on the notation above, the theorem of Sklar is :

\begin{theorem} \label{theo1} For any \textit{cdf} $F$ on $\mathbb{R}^d$, $d\geq 1$, there exists a copula $C$ on 
$\mathbb{R}^d$ such that

\begin{equation}
\forall x \in \mathbb{R}^d, \ F(x)=C(F_1(x),...,F_d(x)). \label{sklar} 
\end{equation}
\end{theorem} 

\bigskip \noindent This theorem is now among the most important tools in Statistics since it allows to study the dependence between the components of a random vector through the copula, meaning that the intrinsic dependence does not depend on the margins.\\

\noindent This note aims to give proof of some lines based on generalized inverses of the margins.\\

\section{A proof of the Sklar's Theorem}

\noindent Define for $s=(s_1,s_2,\cdots,s_d) \in [0,1]^d$,

\begin{equation}
C(s)=F( F_1^{-1}(s_1+0), F_2^{-1}(s_2+0), \cdots, F_d^{-1}(s_d+0)). \label{cop}
\end{equation}
 
\bigskip \noindent It is immediate that $C$ assigns non-negative volumes to cuboids of $[0,1]^d$, since according to Condition (DF2), Formula \eqref{VP} for $C$ derives from the same for $F$ where the arguments are the form $F_i^{-1}(\circ+0)$, $1\leq i\leq d$.\\

\noindent Also $C$ is right-continuous since $F$ is right-continuous as well as each $F_i^{-1}(\circ+0)$, $1\leq i\leq d$. By passing, this explains why we took the right-limits because the $F_i^{-1}(\circ)$'s are left-continuous.\\

\noindent Finally, by combining Formulas \eqref{FF} and \eqref{cop}, we get the conclusion of Sklar in Formula 
\eqref{sklar}. The proof is finished. $\square$\\


\begin{thebibliography}{99}
\bibitem[Billingsley (1968)]{billingsley} Billingsley, P.(1968). \textit{Convergence of Probability measures}. John Wiley, New-York.\\

\bibitem[Lo(2017)]{ips-mestuto-ang} Lo, G. S. (2017) Measure Theory and Integration By and For the Learner.
 \text {SPAS Books Series}. Saint-Louis, Senegal - Calgary, Canada. Doi :  http://dx.doi.org/10.16929/sbs/2016.0005,
ISBN : 978-2-9559183-5-7.\\

\bibitem[Lo \textit{et al.}(2016a)]{ips-wcia-ang} Lo, G.S.(2016). Weak Convergence (IA). Sequences of random vectors. \text {SPAS Books Series}. Saint-Louis, Senegal - Calgary, Canada. Doi : 10.16929/sbs/2016.0001. Arxiv : 1610.05415. ISBN : 978-2-9559183-1-9.\\  

\bibitem[Lo \textit{et al.}(2016b)]{ips-wcia-fr} Lo, G.S.(2016). Convergence vague (IA). Suites de vecteurs al\'eatoires. \text {SPAS Books Series}. Saint-Louis, SENEGAL - CANADA, Canada. Doi : 

\bibitem[Resnick (1987)]{resnick} Resnick, S.I. (1987). \textit{Extreme Values, Regular
Variation and Point Processes}. Springer-Verlag, New-York.\\

\bibitem[Nelsen (2006)]{nelsen} Nelsen, R.B. (2006). \textit{An introduction to copula}. Springer-Verlag, New-York.

\end{thebibliography}
\end{document}